\documentclass[11pt,final,reqno]{amsart}
\usepackage{amsthm}
\usepackage{amsmath}
\usepackage{mathrsfs}
\usepackage{amssymb}
\usepackage{bbm}
\usepackage{enumerate}
\usepackage{array}
\usepackage[pdfstartview=FitH,backref,colorlinks]{hyperref}
\usepackage{enumitem}
\usepackage{graphicx}

\newcommand{\brac}[1]{\left(#1\right)}
\newcommand{\bfrac}[2]{\brac{\frac{#1}{#2}}}
\newcommand{\beq}[1]{\begin{equation}\label{#1}}
\newcommand{\eeq}{\end{equation}}
\newcommand{\blem}[1]{\begin{lemma}\label{#1}}
\newcommand{\elem}{\end{lemma}}
\newcommand{\bth}[1]{\begin{theorem}\label{#1}}
\newcommand{\enth}{\end{theorem}}
\newcommand{\brem}[1]{\begin{remark}\label{#1}}
\newcommand{\erem}{\end{remark}}

\def\eps{\varepsilon}

\newtheorem{theorem}              {Theorem}
\newtheorem{lemma}      [theorem] {Lemma}

\begin{document}

\title{On independent sets in hypergraphs}

\author[]{A.~Kostochka}
\address{Department of Mathematics,
University of Illinois, Urbana-Champaign.}
\email{kostochk@uiuc.edu}
\thanks{A.K.\ was supported by  NSF grant DMS-0965587 and by grant 09-01-00244-a
of the Russian Foundation for Basic Research}

\author[]{D.~Mubayi}
\address{Department of Mathematics, Statistics, and Computer Science, University of Illinois, Chicago, IL 60607}
\thanks{D.M.\ was supported in part by NSF grants DMS 0653946 and DMS 0969092}
\email{mubayi@math.uic.edu}

\author[]{J.~Verstra\"ete}
\address{Department of Mathematics,
University of California, San Diego (UCSD),
La Jolla, CA 92093-0112, USA.}
\email{jacques@ucsd.edu}
\thanks{J.V.\ was supported by NSF grant DMS-0800704}

\thanks{This research originated when the authors participated in the
  AIM workshop \emph{Hypergraph Tur\'{a}n Problem},
  March 21-25, Palo Alto, CA}

\begin{abstract}
The {\em independence number} $\alpha(H)$ of a hypergraph $H$ is the size of a largest set of vertices containing no edge of $H$.
In this paper,  we prove that if $H_n$ is an $n$-vertex $(r + 1)$-uniform hypergraph in which
every $r$-element set is contained in at most $d$ edges, where $0<d < n/(\log n)^{3r^2}$, then
\[ \alpha(H_n) \geq c_r \Bigl(\frac{n}{d}\log \frac{n}{d}\Bigr)^{1/r}\]
  where $c_r > 0$ satisfies $c_r \sim r/e$ as $r \rightarrow \infty$.
The value of $c_r$ improves and generalizes several earlier results.
 Our relatively short proof extends a method due to Shearer.

 The above statement is close to best possible, in the sense that
for each $r \geq 2$ and all values of $d \in \mathbb N$, there are infinitely many $H_n$ such that
\[ \alpha(H_n) \leq b_r \Bigl(\frac{n}{d}\log \frac{n}{d} \Bigr)^{1/r}\]
where $b_r > 0$ depends only on $r$. In addition, for many values of $d$ we show $b_r \sim c_r$ as $r \rightarrow \infty$, so the result is almost sharp for large $r$. We give an application to hypergraph Ramsey numbers involving independent neighborhoods.
\end{abstract}

  \maketitle

\section{Introduction}

In this paper, an {\em $r$-graph} is a set of $r$-element subsets of a finite set, where the sets are called edges and the elements of the finite set are called vertices. An independent set in an $r$-graph is a set of vertices containing no edge. The independence number $\alpha(H)$ of an $r$-graph $H$ is the maximum size of an independent set in $H$.

A {\em partial Steiner $(n,r+1,r)$-system} is an $n$-vertex $(r + 1)$-graph such that each $r$-element set of vertices
is contained in at most one edge. The {\em maximum $r$-degree} of an $(r + 1)$-graph $H$ is the
maximum number of edges that any $r$-set of vertices is contained in.

\medskip

The independence number $\alpha(H)$ has been studied at length in Steiner systems, sometimes in the language of projective geometry,
in terms of {\em maximum complete arcs}, and has applications to geometric problems, for instance the ``orchard planting problem'' (see~\cite{Fu,FP})
or Heilbronn's celebrated triangle problem~\cite{KPS}. Given a partial Steiner $(n,r+1,r)$-system $H$, Phelps and R\"{o}dl~\cite{PR} were the first to show $\alpha(H) > c(n\log n)^{1/r}$ for some constant $c > 0$ depending only on $r$, answering a question of Erd\H{o}s~\cite{E}. R\"odl and \v{S}inajov\'a~\cite{RS} proved that this result is tight apart from the constant $c$

One of the methods for finding large independent sets is the randomized greedy approach: one picks a small set of independent vertices repeatedly, delete
the neighbors of this set, and control the statistics of the remaining hypergraph at each stage. The paper of Ajtai, Koml\'{o}s, Pintz, Spencer and Szemer\'{e}di~\cite{AKPSS} gives a detailed analysis of such an algorithm for finding independent sets in $r$-graphs with bounded degrees.
This approach has been used successfully to attack the corresponding coloring problems for hypergraphs (see~\cite{BFM,FM, FM2}).

\subsection{Main Theorem}

 In this paper, we give a short proof of a general result for $(r + 1)$-graphs with maximum $r$-degree $d$.
  This extends the afore-mentioned result of Phelps and R\"{o}dl, which is the case $d = 1$, without a randomized greedy approach.
 Shearer gave an ingenious short proof of the celebrated result of Ajtai, Koml\'os, Szemer\'edi~\cite{AKS} that every triangle-free graph with $n$ vertices and average degree $d$ has an independent set of size at least $c (n/d) \log d$ for some constant $c$.  He asked whether his method could be applied to the hypergraph setting and we partially answer his question by proving our main result using his approach:

\begin{theorem}\label{main}
Fix $r \ge 2$. There exists $c_r > 0$ such that if $H$ is an $(r + 1)$-graph on $n$ vertices with maximum $r$-degree $d < n/(\log n)^{3r^2}$, then
\[ \alpha(H) \geq c_r \Bigl(\frac{n}{d}\log \frac{n}{d}\Bigr)^{\frac{1}{r}}\]
where $c_r > 0$ and $c_r \sim r/e$ as $r \rightarrow \infty$.
\end{theorem}

Theorem \ref{main} is close to best possible
as for any values of $r \geq 2$, there exists an $(r + 1)$-graph $H$ on $n$ vertices with maximum $r$-degree $d$ and, for some constant $b_r$,
\[ \alpha(H) \leq b_r \Bigl(\frac{n}{d}\log \frac{n}{d}\Bigr)^{\frac{1}{r}}.\]
Furthermore, if $d \gg \log n$ and $\log d \ll \log n$, then we show $b_r \sim r/e \sim c_r$ as $r \rightarrow \infty$, in Section 3, so in this range of $d$ and for $r \rightarrow \infty$, Theorem \ref{main} is best possible including the constant. The best constant $c_r$ that can be read out of the proof of Theorem \ref{main} is
\[ c_r = \Bigl(\frac{r!}{r(3r - 1)2^r\log(1 - 2^{-r})}\Bigr)^{1/r}\]
and gives $c_3 \approx 0.538$. This is the current best lower bound on the independence number of a Steiner triple system.
An upper bound of $4\sqrt{n\log n}$ for Steiner triple systems was given by Phelps and R\"{o}dl~\cite{PR} and generalized to
Steiner $(n,r,k)$-systems by R\"{o}dl and \v{S}inajov\'{a}~\cite{RS}.

\medskip

\subsection{Independent neighborhoods}

An $r$-graph $H$ is said to have {\em independent neighborhoods} if for every set $R$ of $r - 1$ vertices,
$\{e \backslash R : R \subset e \in H\}$ is an independent set. These hypergraphs have been studied from
the point of view of extremal hypergraph theory~\cite{FMP,FPS} and hypergraph coloring~\cite{BFM}.
Denote by $T_r$ the $r$-graph with vertex set $R \cup S$ with
$|R| = r$ and $|S| = r - 1$ and consisting of all edges containing $S$ together with the edge $R$. Then an $r$-graph has
independent neighborhoods if and only if it does not contain $T_r$ as a subgraph. The {\em Ramsey number} $R(T_r,K_t^{(r)})$ is the minimum
$N$ such that in every red-blue coloring of the edges of the complete $r$-graph $K_N^{(r)}$ on $N$ vertices, there is either a red $T_r$ or a blue $K_t^{(r)}$.
As a straightforward consequence of Theorem \ref{main}, we obtain the following result:

\begin{theorem}\label{idn}
Let $H$ be an $r$-graph on $n$ vertices with independent neighborhoods. Then for some constant $c$,
$\alpha(H) \geq c(n\log n)^{\frac{1}{r}}$. In particular,
\[ R(T_r,K_t^{(r)}) = O\Bigl(\frac{t^r}{\log t}\Bigr).\]

\end{theorem}
We remark here that the bound above without the log factor is trivial and it follows from known results that $R(T_r, K_t^{(r)})> c t^{r}/(\log t)^{r/(r-1)}$ for suitable $c>0$ depending on $t$.
We believe that the Ramsey result is best possible up to the value of the implicit constant.
In the case $r = 2$, for graphs, a graph has independent neighborhoods if and only if it is triangle-free. Theorem
\ref{idn} therefore generalizes the well-known result of Ajtai, Koml\'{o}s and Szemer\'{e}di~\cite{AKS} for triangle-free graphs to hypergraphs.
It remains an open problem to show that Theorem \ref{idn} is  best possible for all $r$. It is known to be best possible for graphs by a result of Kim~\cite{K} which establishes that $R(K_3,K_t^{(r)})$ has order of magnitude $t^2/(\log t)$.

\medskip

\subsection{Organization}

This paper is organized as follows: we start with stating the Chernoff Bound in Section 2, which will be used repeatedly in the probabilistic methods
to follow. In Section 3, we give the constructions which prove that Theorem \ref{main} is tight up to the constant $c_r$. In Section 4 we will sketch the proof for the case $r=2$ -- the interested reader might want to read this section first to see the main ideas.
In Sections 5 and 6,
we establish some preliminaries for the proof of the Theorem \ref{main}, which is in Section 7.
In Section 8 we give an application to Ramsey numbers and hypergraphs with independent neighborhoods.
We end with some concluding remarks.

\medskip

\subsection{Notation}

A hypergraph $H$ is a pair $(V(H), E(H))$ where $E(H) \subset 2^{V(H)}$; it is an $r$-graph if $E(H) \subset {V(H) \choose r}$.  Sometimes we will abuse notation by associating $H$ with its edge set $E(H)$.  A {\em triangle} in an $r$-graph $H$ is a subgraph of three edges $\{e,f,g\}$ such that
$|e \cap f| = |f \cap g| = |g \cap e| = 1$ and all the intersections are distinct.
 A hypergraph is {\em linear} if it has no pair of distinct edges sharing two or more vertices.
 A set $Z \subseteq V$ is an {\em independent set} of $H$ if $Z$ contains no edges of $H$.
Two vertices of $H$ are {\em adjacent} if they are contained in a common edge of $H$. Let $N(x)$ denote the
set of vertices adjacent to $x \in V(H)$. A {\em subgraph} or {\em subhypergraph} of a hypergraph $H = (V,E)$
is a hypergraph $H' = (V',E')$ where $V' \subseteq V$ and $E' \subseteq E$.
For $X \subset V$, the subgraph of $H$ {\em induced} by $X$ is the subgraph $H[X]$ consisting
of all edges of $H$ that are contained in $X$.

All logarithms in this paper are to the natural base, $e$.
We write $f(n) \sim g(n)$  or $f(n)=(1+o(1))g(n)$ for functions $f,g : \mathbb N \rightarrow \mathbb R^+$ to denote $f(n)/g(n) \rightarrow 1$ as $n \rightarrow
\infty$, and $f(n) = O(g(n))$ to denote that there is a constant $c$ such that $f(n) \leq cg(n)$ for all $n$.
 We also write $f(n) \lesssim g(n)$ if $\limsup f(n)/g(n) \leq 1$ as $n \rightarrow \infty$.  Similarly, $f(n) \ll g(n)$ or $f(n)=o(g(n))$ means that $\lim f(n)/g(n)=0$.
Unless otherwise indicated, any asymptotic notation implicitly assumes $n \rightarrow \infty$.

\medskip

\section{Chernoff-type bounds}

The proof of Theorem \ref{main} is probabilistic. In the subsequent material, we shall make use of the following concentration inequality, which is a generalization of the Chernoff Bound (see McDiarmid: Theorem 2.7 in~\cite{McDiarmid}).
In this section, $U \sim \mbox{binomial}(n,p)$ means $U$ is a binomial random variable with success probability $p$ in $n$ trials. Throughout the paper, if $(A_n)_{n \in \mathbb N}$ is a sequence of events in some probability space, then we say $A_n$ occurs {\em with high probability} if $\lim_{n \rightarrow \infty} P(A_n) = 1$.

\begin{lemma}\label{C}
Let $U$ be a sum of independent random variables $U_1,U_2,\dots,U_n$ such that $E(U) = \mu$ and
$U_i \leq E(U_i) + b$ for all $i$. Let $V$ be the variance of $U$. Then for any $\lambda>0$
\begin{itemize}
\item[1.] $P(U  \geq \mu+ \lambda) \leq e^{-\frac{\lambda^2}{2V+ b\lambda}}$.
\item[2.] If $U \sim \mbox{binomial}(n,p)$, then
$P(|U - \mu| \geq \varepsilon \mu) \leq 2e^{-\frac{\varepsilon^2 \mu}{2}}.$
\end{itemize}
\end{lemma}

The inequality in Lemma \ref{C} part 2 will be referred to as the Chernoff Bound~\cite{Cher}.

\medskip

\subsection{A technical lemma} In the proof of Theorem \ref{main}, we require the following consequence
of the Chernoff Bound:

\begin{lemma}\label{consequence}
Let $k,b$ be positive integers and $q \in (0,1]$, and define
\begin{equation}\label{S}
S := \sum_{j = 0}^k {k \choose j} q^j(1-q)^{k-j} \min\{j,b\}.
\end{equation}
Then as $k \rightarrow \infty$,
\begin{equation}\label{Sbound}
S \sim \min\{qk,b\}.
\end{equation}
\end{lemma}

\begin{proof}
Let $Y \sim \mbox{Bin}(k, q)$. Then clearly
\begin{eqnarray*}
S = \sum_{j = 0}^k P(Y = j) \min\{j,b\} = E(\min\{Y,b\}).
\end{eqnarray*}
By Lemma \ref{C} part 2, $Y \sim qk$ with high probability as $k \rightarrow \infty$. Therefore $E(\min\{Y,b\})= (1-o(1))\min\{qk, b\}+o(1)b \sim \min\{qk,b\}$.
\end{proof}

\section{Hypergraphs with low independence numbers}

We show that Theorem \ref{main} is tight for all $d \in \mathbb N$ up to the value of the constant $c_r$, using a ``blowup'' of a Steiner system. Furthermore,
for many values of $d$ and large $r$, we shall see via a random hypergraph construction that the constant $c_r$ is itself almost best possible.

\subsection{Blowup of a Steiner system}
 Let $S_n$ be any Steiner $(n,r+1,r)$-system
with $V(S_n) = \{1,2,\dots,n\}$. Define a hypergraph $H=(V,E)$ with $N = dn$ vertices and
with maximum $r$-degree $d$ as follows: let $V$ be a disjoint union of sets $V_1,V_2,\dots,V_n$ each of size $d$. For each edge $e=\{x_1, \ldots, x_r\} \in S_n$ let $B_e$ be the collection of all edges of the form $\{v_1, \ldots, v_r\}$ where $v_i \in V_{x_i}$.
Let $E$ comprise
all edges in each $V_i$ together with all edges in each $B_e$. Note that
every edge $e \in H$ has the property that either $e \subset V_i$ for some $i$ or $|e \cap V_i| = 1$ for exactly $r$ values of $i$. We may refer to $H$ loosely as a blowup of a Steiner system.
We observe that $\alpha(H) = r\alpha(S_n)$ since every independent set $X$ of $H$ contains at most $r$ vertices in each $V_i$, and $\{i : |X \cap V_i| \neq \emptyset\}$ is an independent set
of $S_n$. It is known that there are Steiner $(n,r+1,r)$-systems $S_n$ in which $\alpha(S_n) \leq a_r(n\log n)^{1/r}$ for some $a_r > 0$ depending only on $r$ -- see~\cite{PR,RS}. Therefore blowing up these Steiner systems, we obtain
$(r + 1)$-graphs $H$ with $N$ vertices and maximum $r$-degree $d$ such that
\begin{eqnarray*}
\alpha(H) &=& r\alpha(S_n) \\
&\leq& ra_r \Bigl(\frac{N}{d}\log \frac{N}{d}\Bigr)^{1/r} \\
&\leq& b_r \Bigl(\frac{N}{d}\log \frac{N}{d}\Bigr)
\end{eqnarray*}
where $b_r > 0$ depends only on $r$. This shows Theorem \ref{main} is tight up to the constant $c_r$.

\subsection{Random hypergraphs}

A natural candidate for an $(r + 1)$-graph with low independence number is the random $(r + 1)$-graph $H = H_{n,r+1,p}$.
This probability space is defined by selecting randomly and independently
with probability $p$ edges of the complete $r$-uniform  hypergraph on $n$ vertices, and letting $H$ be the $(r + 1)$-graph of selected edges. We sketch a standard argument showing that a random hypergraph
gives good examples of a hypergraph with low independence number. We take $p = d/(n - r)$, so that the expected $r$-degree of any $r$-element set in $V(H)$ is exactly $d$. By the Chernoff Bound, Lemma \ref{C}.2,
if $d \gg \log n$, then with high probability, every $r$-set in $H$ has $r$-degree asymptotic to $d$. Next,
using the bounds $(1 - p)^y \leq e^{-py}$ for $p \in [0,1]$ and $y \geq 0$ and $(a-b + 1)^b/b! \leq {a \choose b} \leq a^b$ for $a \geq b \geq 1$,
the expected number of independent sets of size $x$ in $H$ is exactly
\[ E := {n \choose x} (1 - p)^{{x \choose r + 1}} < \exp\Bigl(x\log n - \frac{d}{n} \cdot \frac{(x - r)^{r + 1}}{(r + 1)!}\Bigr).\]
Fix $\eps>0$ and let
$$x=(1+\eps)(r+1)!^{1/r}\left(\frac{n}{d}\log n\right)^{1/r}.$$
Then, as $n \rightarrow \infty$, we see that
$$\frac{d}{n} \frac{ (x-r)^{r+1}}{(r+1)!}>x \log n$$
and therefore $E<1$.
We conclude that with positive probability, $\alpha(H) < x$ and consequently,
$$\alpha(H) \lesssim (r + 1)!^{1/r} \Bigl(\frac{n}{d}\log n\Bigr)^{1/r}$$
as required.
If, in addition, $\log d \ll \log n$, then $\log \frac{n}{d} \sim \log n$ and so
$$\alpha(H) \lesssim (r + 1)!^{1/r} \Bigl(\frac{n}{d}\log \frac{n}{d}\Bigr)^{1/r}.$$
Note that $(r + 1)!^{1/r} \sim r/e \sim c_r$ showing that Theorem \ref{main} provides close to the right constant for large $r$.

\medskip

\section{Sketch Proof of Theorem \ref{main}}

We outline the proof of Theorem \ref{main} for linear triple systems -- that is when $r = 2$ and $d = 1$ -- since
the general proof requires only slight modifications of the ideas in this case.
For a contradiction, suppose there are $n$-vertex linear triple systems $H$ such that $\alpha(H) \ll \sqrt{n\log n}$.

\subsection{Step 1 : Random sets} A {\em random set} is a set $X \subset V(H)$ whose vertices are chosen independently from
$H$ with probability $$p = \frac{n^{-2/5}}{(\log\log \log n)^{3/5}}.$$ Then $E(|X|) = pn$ and $E(|T|) \leq p^6{n \choose 3}$ where $T=T(X)$ is the set of triangles in $H[X]$.  The second bound holds since a triangle is uniquely determined by the three vertices which are the pairwise intersections of its edges, since $H$ is linear. The choice of $p$ ensures $E(|T|) \ll pn$. For an independent set $Z \subset V(H)$ and $x \in X$, let
$$\omega_Z(x)= \min(\, \log n, \, |\{xyz \in E(H): \{y,z\} \subset Z\}|\, ).$$
Define
\[ h(Z,X) = \sum_{x \in X \backslash Z} \omega_Z(x).\]
Since $H$ is linear, each $\{y,z\} \subset Z$ accounts for at most one such triple $\{x,y,z\}$
and $x \in X$ with probability $p$, so $$E(h(Z,X)) \leq p{|Z| \choose 2} \leq p \alpha(H)^2 \ll pn\log n.$$
We use Lemma \ref{C} -- details are given in Section 6 -- to show that $X$ can be chosen so that \medskip

1) $h(Z,X) \ll pn\log n$ for all independent sets $Z$ in $H$,

 2) $|X| \sim pn$
and

3) $H[X]$ is linear and $T(X) = 0$.
\medskip

 \noindent
 Henceforth, fix such a subset $X$ and work in $H[X]$.

\subsection{Step 2 : Random weights} Let $Z$ be a {\em randomly and uniformly chosen independent set} in $H[X]$ and define
for $x \in X$ the random variable
\[ W_x = \left\{\begin{array}{ll}
p\sqrt{n} & \mbox{ if }x \in Z\\
\omega_Z(x) & \mbox{ if }x \in X \backslash Z
\end{array}\right.\]
We bound the expected value of $W := \sum_{x \in X} W_x$ in two ways.

\subsection{Step 3 : Upper bound for random weights} By definition we have
$W = p\sqrt{n}|Z| + h(Z,X)$. The choice of $X$ in Step 1 ensures that $$W \le  p\sqrt{n}\alpha(H) + o(pn\log n)=o(pn\log n)$$
 so $E(W) \ll pn\log n$.

\subsection{Step 4 : Lower bound for random weights} Fixing an $x \in X$, we condition on the value of $Z_x = Z \backslash (N(x) \cup \{x\})$.
Fixing $Z_x$, let $J \subset N(x)$ be the set of vertices such that $Z_x \cup J$ is an independent set in $H[X]$.
Since $H[X]$ is triangle-free and linear, no edge of $H[X]$ has two vertices in $N(x)$ except the edges on $x$.
Therefore, for any independent set $I$ in $H[J \cup \{x\}]$, $I \cup Z_x$ is an independent set. Let $M$
 be the set of pairs of vertices of $J$ forming an edge with $x$ and $L$ be the set of vertices in $J$ not incident to any pair of $M$.
  If  $|M| = k$, then there are
$4^k + 3^k$ independent sets in $H[\bigcup M \cup \{x\}]$ -- those not containing $x$ plus those containing $x$ --
and by the definition of $W_x$
\[ E(W_x|Z_x) = \frac{2^{|L|}p\sqrt{n} 3^k + 2^{|L|}\sum_{j = 0}^k {k \choose j} 3^{k - j} \min\{j,\log n\}}{2^{|L|}(3^k + 4^k)}.\]
Using Lemma \ref{consequence}, with $q = 1/4$, the sum is asymptotic to $\min\{k4^{k - 1}, 4^k \log n\}$ if $k \rightarrow \infty$.
By the choice of $p$, a calculation shows the minimum value of the right hand side is of order
$\log n$ -- see Section 5 for details. So for every $x \in X$, $E(W_x|Z_x) = \Omega(\log n)$. Therefore by the tower property,
\[ E(W) = \sum_{x \in X} E(W_x) = \sum_{x \in X} E(E(W_x|Z_x)) = \Omega(pn\log n).\]
This contradicts the upper bound in Step 3, and completes the proof.

\section{An inequality on independent sets}

It will be shown that if $H$ is an $(r + 1)$-graph of maximum $r$-degree $d$,
then $H$ has a large linear triangle-free subgraph, and that subgraph will contain an independent set
of the size stated in Theorem \ref{main}. In this section, we prove a general inequality for independent sets in linear-triangle-free
$r$-graphs. Let $H$ be a linear triangle-free $(r + 1)$-graph with $m$ vertices.
Let $\mathcal{Z}$ be the set of all independent sets of $H$. The key quantity
we wish to control is defined as follows. For $Z \in \mathcal{Z}$, $b \in \mathbb R$, and $v \in V(H) \backslash Z$, define $\omega_Z(v,b)$ to be the minimum of $b$ and the number of $r$-sets $e \subset Z$ such that $e \cup \{v\} \in H$. Then define
\[ h(Z,b) = \sum_{v \in V(H) \backslash Z} \omega_Z(v,b).\]

\begin{lemma}\label{inequality}
Let $H$ be a linear triangle-free $(r + 1)$-graph with $m$ vertices, and let $Z$ be a uniformly randomly chosen independent set
in $H$, and $b \in \mathbb R^+$. Then as $b \rightarrow \infty$,
\begin{equation}\label{maininequality}
E(h(Z,b)) + e^b E(|Z|) \gtrsim \frac{bm}{-2^r\log(1 - 2^{-r})}.
\end{equation}
\end{lemma}

\medskip

\begin{proof}
Let $V = V(H)$ and $q = 1 - 2^{-r}$. For $v \in V$, define the random variable:
\[ W_v = \left\{\begin{array}{ll}
e^b & \mbox{ if }v \in Z \\
\omega_Z(v,b) & \mbox{ if }v \in V \backslash Z
\end{array}\right.
\]
By definition of $W_v$,
\[ W := \sum_{v \in V(H)} W_v = \sum_{v \in Z} W_v + \sum_{v \in V \backslash Z} W_v = e^b |Z| + h(Z,b).\]
To complete the proof, we show $E(W_v) \gtrsim b/(-2^r\log q)$ for every $v \in V$.

\medskip

Fixing $v \in V$ and $Z_v = Z \backslash (N(v) \cup \{v\})$, define
\[ J = \{u \in N(v) : Z_v \cup \{u\} \in \mathcal{Z}\}.\]
Since $H$ is linear and triangle-free, $Z$ is obtained from $Z_v$ by selecting an independent subset of $H[J \cup \{v\}]$.
Let $M$ be the set of $r$-sets in $J$ forming an edge with $v$ and let $L=J-\bigcup M$.
Since $H$ is linear, $M$ consists of disjoint $r$-sets. A set of vertices of $J \cup \{v\}$ containing $v$ is independent in $H$ if and only if it contains at most $r - 1$ vertices from each of the sets in $M$ together with any subset of $L$. Any independent set of $H$ in $J \cup \{v\}$
not containing $v$ consists of any subset of $\bigcup M \cup L$. If $|M| = k$ and $|L| = \ell$, there are $2^{\ell}(2^{rk} + (2^r - 1)^{k})$ independent sets in $H[J \cup \{v\}]$. It follows from the definition of $W_v$ that
\begin{eqnarray}
 E(W_v|Z_v) &=& \frac{e^b2^{\ell}(2^r - 1)^k + 2^{\ell}\sum_{j = 0}^k {k \choose j} (2^r - 1)^{k - j}\min\{j,b\}}{2^{\ell}(2^{rk} + (2^r - 1)^k)} \nonumber \\
&=& \frac{e^b q^k}{1 + q^k}+ \frac{\sum_{j = 0}^k {k \choose j} (2^r - 1)^{k - j}\min\{j,b\}}{2^{rk} + (2^r - 1)^k)}. \label{expectation}
\end{eqnarray}
We shall show $E(W_v|Z_v) \gtrsim b/(-2^r\log q)$. First suppose that $e^b q^k > 2b$. Then using the inequality $-\log(1-x)>x$ for $0<x<1$, we obtain
\[
E(W_v|Z_v) \geq \frac{e^b q^k}{1 + q^k} >\frac{e^bq^k}{2}> b > \frac{b}{-2^r \log q}.
\]
Next suppose that $e^b q^k \leq 2b$. Then
Lemma \ref{consequence} gives
$$\sum_{j = 0}^k {k \choose j} (2^r - 1)^{k - j}\min\{j,b\} \sim 2^{rk} \min\{(1-q)k, b \}.$$
Consequently,
 $$
 E(W_v|Z_v) \gtrsim \frac{ 2^{rk} \min\{(1-q)k, b \}}{2^{rk} + (2^r - 1)^k}.$$
Since $k \rightarrow \infty$ as $b \rightarrow \infty$,

$$\frac{ 2^{rk} \min\{
(1-q)k, b \}}
{2^{rk} + (2^r - 1)^k} \sim
\min\{(1-q)k, b \}$$
Since $e^b q^k \leq 2b$, we have $k>(\log 2b-b)/\log q \sim  -b/\log q$, and so
$$\min\{(1-q)k, b\} \gtrsim \min\left \{ \frac{(1-q)b}{-\log q}, b\right\} = \min\left\{ \frac{b}{-2^r \log q}, b \right\} \ge \frac{b}{-2^r \log q}$$
Now (\ref{expectation}) and the tower property of expectation implies,
\[
E(W) = \sum_{v \in V} E(E(W_v|Z_v)) \gtrsim \frac{bm}{-2^{r}\log q}.
\]
This completes the proof of Lemma \ref{inequality}.
\end{proof}

\bigskip

\section{Random subsets of hypergraphs}

To prove Theorem \ref{main}, we shall find an appropriate
set $Y \subset V(H)$ such that $H[Y]$ is linear and triangle-free and then we apply Lemma \ref{inequality}.
To do so, we need to find a set $Y$ in which the quantity $h(Z, b)$ in Lemma \ref{inequality} is not too large.
 The set $Y$ will be found by random sampling. A {\em random set} refers to a set $X \subset V(H)$ whose vertices are chosen from $V$ independently with probability $p$, where $p$ is to be chosen later.

\begin{lemma}\label{goodX}
Let $H$ be an $n$-vertex $(r + 1)$-graph with maximum $r$-degree $d$ and $\alpha(H) \leq \alpha$. Suppose that for some $p \in [0,1]$ with $p\gg 1/n$ and $b \in \mathbb R^+$,
\begin{equation}\label{condition}
\frac{pd^2\alpha^{2r}}{nb^2 + db\alpha^r} \gg \alpha \log n \quad \mbox{ and } \quad d^3 n^{3r - 3}p^{3r} \ll pn.
\end{equation}

Then there exists a set $Y \subseteq V(H)$ with the following properties

 $\bullet$ $|Y| \sim pn$

  $\bullet$ $H[Y]$ is linear and triangle-free and

  $\bullet$ for every independent
set $Z$ in $H[Y]$,
\begin{equation}\label{zbound}
h(Z, b)\lesssim pd {\alpha \choose r}.
\end{equation}
\end{lemma}

\begin{proof}
Let $X$ be a random subset of $V:=V(H)$. The main part of the proof is to show that with high probability, $h(Z, b) \lesssim pd{\alpha \choose r}$ for every independent set $Z$ in $H[X]$. First we upper bound $E(h(Z, b))$. Since $H$ has maximum $r$-degree $d$,
\[ E(h(Z, b)) \leq dp{|Z| \choose r} \leq dp{\alpha \choose r}\]
for any independent set $Z$ in $H$. Now $h(Z, b)$ is a sum of independent random variables $\omega_Z(v,b)$,
each bounded by $b$. Letting $I_{v \in X}$ be the indicator that $v$ is in $X$, we have:

\begin{eqnarray*}
Var(h(Z,b)) &\leq& E(h(Z,b)^2) \\
&\leq& \sum_{v \in V \backslash Z} E(I_{v \in X}\omega_Z(v,b)^2) \\
&\leq& \sum_{v \in V} E(I_{v \in X})b^2 \\
&=& pnb^2.
\end{eqnarray*}

By Lemma \ref{C} part 1 with
$\varepsilon > 0$, $\lambda = \varepsilon pd{\alpha \choose r}$,
\begin{eqnarray*}
-\log P(h(Z,b) > E(h(Z,b)) + \lambda) &\geq& \frac{\lambda^2}{2pnb^2 + \lambda b} \\
&=& \frac{(\varepsilon pd)^2 {\alpha \choose r}^2}{2pnb^2 + \varepsilon pd {\alpha \choose r} b} \\
&\geq& \frac{(\varepsilon pd\alpha^r)^2}{3r!^2(pnb^2 + pdb\alpha^r)} \; \; \gg \; \; \alpha \log n
\end{eqnarray*}
by (\ref{condition}). Since $|\mathcal{Z}| < n^{\alpha(H)}$, this
shows by Markov's Inequality that with high probability, $h(Z,b) \leq (1 + \varepsilon)pd{\alpha \choose r}$.
Since $\varepsilon > 0$ is arbitrary, this means $h(Z,b) \lesssim pd{\alpha \choose r}$.

\medskip

Consider pairs of edges in $X$ that intersect in at least two vertices.  The number of pairs of edges in $H$ that intersect in $i$ vertices can be upper bounded as follows: First choose an $i$-set $S$ of vertices that is the intersection of two edges -- there are at most $n^i$ ways of choosing $S$. Now consider the $(r+1-i)$-graph $H_S$ consisting of edges of the form $E-S$ where $E \in E(H)$.  Since $H$ has $r$-degree at most $d$, we conclude that $H_S$ has $(r-i)$-degree at most $d$, so $H_S$ has at most $dn^{r-i}$ edges. Now we pick two edges in $H_S$ that are disjoint. The number of ways of doing this is at most $d^2 n^{2r-2i}$.  Altogether, the number of pairs of edges in $H$ sharing exactly $i$ vertices is at most $d^2 n^{2r-i}$, and the probability that one such pair lies in $X$ is $p^{2r+2-i}$.
We conclude, using $p \gg 1/n$ and (\ref{condition}), that the expected number of pairs of edges in $X$ intersecting in two or more vertices is at most
\[ d^2(p^{2r}n^{2r - 2} + p^{2r - 1}n^{2r - 3} + \dots + p^{r + 3}n^{r + 1}+p^{r+2}n^r)\ll d^2p^{2r}n^{2r-2} \ll pn.\]

 Next we consider triangles in $H[X]$ which here are triples $\{e,f,g\}$ of edges of $H[X]$
such that $|e \cap f| = |f \cap g| = |g \cap e| = 1$ and $e \cap f \cap g = \emptyset$.  There are fewer than $n^3$ choices
for $e \cap f, f \cap g,g \cap e$. Fixing $e \cap f$ and $e \cap g$, there are fewer than $n^{r - 2}d$ choices for $e$ since $H$ has
$r$-degree at most $d$. It follows that the expected number of triangles in $H[X]$ is $d^3 n^{3r - 3}p^{3r} \ll pn$, using (\ref{condition}).
We conclude that the number of triangles $T=T(X)$ in $H[X]$ satisfies
$E(|T(X)|) \ll pn$. Now if $Y$ is obtained from $X$ by deleting a vertex of $X$ from each triangle in $H[X]$ and from each pair of edges of $H[X]$ intersecting in at least two vertices in $H[X]$, then $|Y| \sim pn$ with high probability. Finally, we observe that the value of $h(Z,b)$ does not increase by deleting vertices from $X$, so (\ref{zbound}) holds in $H[Y]$ with high probability.
\end{proof}

\medskip

\section{Proof of Theorem \ref{main}}

We are now ready to prove Theorem \ref{main}, using Lemmas \ref{inequality} and \ref{goodX}. In the proof, all asymptotic notation refers to $n \rightarrow \infty$. Let $H$ be an $(r + 1)$-graph of maximum $r$-degree $d \leq n/(\log n)^{3r^2}$ on $n$ vertices and independence number at most
\[ \alpha := c\Bigl(\frac{n}{d}\log \frac{n}{d}\Bigr)^{1/r}\]
where $c > 0$ is a constant depending only on $r$. To complete the proof we show $c \geq c_r$ if $n$ is large enough, where
\begin{equation}\label{cr}
c_r^r = \frac{r!}{-r(3r - 1)2^{r}\log(1 - 2^{-r})}.
\end{equation}
This implies that every such $r$-graph has large independence number.
Define $p \in [0,1]$ and $b \in \mathbb R^+$ by
\[ pn = \Bigl(\frac{n}{d\log\log\log n}\Bigr)^{\frac{3}{3r - 1}} \quad \mbox{ and } \quad b = \frac{1}{r(3r - 1)}\log \frac{n}{d}.\]
There are two steps to the proof: first we have to verify that the above choice of parameters allows us to apply Lemma \ref{inequality} and
Lemma \ref{goodX}, in particular (\ref{condition}).

\medskip

We claim that the following hold, which allow us to apply the lemmas:
\begin{eqnarray}
(pd^2 \alpha^{2r})/(nb^2 + db\alpha^r) &\gg& \alpha \log n \label{first} \\
d^3 n^{3r - 3}p^{3r} &\ll& pn \label{second} \\
e^b \alpha &\ll& pd\alpha^r. \label{third}
\end{eqnarray}
The inequality (\ref{second}) follows immediately from the definition of $pn$, due to the $\log\log\log n$ term there. To prove (\ref{first}),
note $nb^2 < db\alpha^r$ and then
\begin{eqnarray*}
\frac{pd^2 \alpha^{2r}}{nb^2 + db\alpha^r} \; \; > \; \; \frac{pd^2\alpha^{2r}}{2db\alpha^r} &=& \frac{pd\alpha^r}{2b} \\
&=& \frac{r(3r - 1) c^r}{2} pn.
\end{eqnarray*}
By the definition of $pn$ and $d \leq n/(\log n)^{3r^2}$, a short calculation yields $pn \gg \alpha \log n$,
which proves (\ref{first}).
For (\ref{third}), we have
\begin{eqnarray*}
e^b \alpha &=& c\Bigl(\frac{n}{d}\Bigr)^{1/r(3r - 1)} \cdot \Bigl(\frac{n}{d}\log\frac{n}{d}\Bigr)^{1/r} \\
&=& c\Bigl(\frac{n}{d}\Bigr)^{3/(3r - 1)} \Bigl(\log\frac{n}{d}\Bigr)^{1/r} \\
&=& c(\log\log\log n)^{3/(3r - 1)} pn \Bigl(\log\frac{n}{d}\Bigr)^{1/r} \; \; \ll \; \; pd\alpha^r
\end{eqnarray*}
since $d \leq n/(\log n)^{3r^2}$ and $r \ge 2$. This verifies (\ref{third}) and so we now apply Lemma \ref{goodX}.

\medskip

By Lemma \ref{goodX}, there is a linear triangle-free subgraph $H[Y]$ with $|Y| \sim pn$ and
\[ h(Z, b) \lesssim pd{\alpha \choose r}\]
for every independent set $Z$ in $H[Y]$. In particular, using (\ref{third}),
\begin{equation}\label{upper}
E(h(Z,b)) + e^b E(|Z|) \lesssim pd{\alpha \choose r} + e^b \alpha \lesssim \frac{c^r}{r!} pn \Bigl(\log \frac{n}{d}\Bigr).
\end{equation}
We note that $b \rightarrow \infty$ since $d \leq n/(\log n)^{3r^2}$. Therefore by Lemma \ref{inequality},
\begin{equation}\label{lower}
E(h(Z,b)) + e^b E(|Z|) \gtrsim \frac{pnb}{-2^{r}\log(1 - 2^{-r})} \gtrsim \frac{c_r^r}{r!} pn \Bigl(\log \frac{n}{d}\Bigr).
\end{equation}
Comparing (\ref{lower}) with (\ref{upper}) gives $c \gtrsim c_r$, as required.

\section{Ramsey numbers and independent neighborhoods}

A celebrated paper of Kim~\cite{K} together with an earlier upper bound of Ajtai, Koml\'{o}s and Szemer\'{e}di~\cite{AKS} shows
that the Ramsey number $R(3,t)$ has order of magnitude $t^2/(\log t)$. Using Theorem \ref{main}, we can generalize part of this result to
hypergraphs in the following manner. Let $T_r$ denote the $r$-graph consisting of $r$ edges containing a given $(r - 1)$-element set,
together with one further edge disjoint from that set and containing one vertex from each of the $r$-edges.
Theorem \ref{idn} is an easy consequence of Theorem \ref{main}:

\begin{proof} [Proof of Theorem \ref{idn}.]
Let $t$ be the lower  bound on $\alpha(H)$. If $H$ has maximum $(r - 1)$-degree at least $t$, then
the set of vertices adjacent to an $r$-set of degree $t$ is an independent set, since $H$ has independent neighborhoods.
Otherwise, by Theorem \ref{main},
\[ \alpha(H) \geq c\Bigl(\frac{n\log n}{t}\Bigr)^{\frac{1}{r - 1}}\]
for an appropriate constant $c$.
  A short computation shows this gives the required upper bound on Ramsey numbers.
\end{proof}

The above theorem is best possible for $r = 2$, as shown via a random construction of triangle-free graphs by Kim~\cite{K}. We believe Theorem \ref{idn} is best
possible for $r > 2$ as well. It is straightforward to give an example with $\alpha(H) \leq c'n^{1/r}(\log n)^{1/(r - 1)}$ with $c' > 0$ using the random hypergraph $H_{n,p}$ with edge probability
$p \approx n^{-(r - 1)/r}$.  One can then use the Local Lemma or the deletion method (see the proof of Theorem 4 in \cite{BFM} for details using the latter approach).

\bigskip

\section{Concluding remarks}

$\bullet$  Duke, Lefmann and R\"{o}dl~\cite{DLR}, based on a paper of Ajtai, Koml\'{o}s, Pintz, Spencer and Szemer\'{e}di~\cite{AKPSS} showed that a linear $(r + 1)$-graph on $n$ vertices with averaged degree $d$ has an independent set of size at least $c'n(\frac{\log d}{d})^{1/r}$. It would be interesting to find a way to
extend the method of this paper to prove such a result.

\medskip

$\bullet$ This paper was partly inspired by the conjecture by Frieze and Mubayi~\cite{FM} that if $H$ is an $(r + 1)$-graph on $n$ vertices with maximum degree $d$,
and $H$ does not contain a specific $(r + 1)$-graph $F$, then $H$ has chromatic number $O(d^{1/r}/(\log d)^{1/r})$. It is generally thought
that proving upper bounds on the chromatic number under such restrictions is a much more difficult problem than finding a large independent set.
In~\cite{FM} and \cite{FM2} the case when $H$ is linear is dealt with using a randomized greedy approach. For $r = 1$ -- i.e. for graphs -- this is known to be true when
$F$ is a bipartite graph, or one vertex away from a bipartite graph~\cite{AKSu}. It is open for graphs even in the case $F = K_4$, and in each case
where the chromatic number conjecture is open, the corresponding question for independence number is also open.

\section{Acknowledgments}
We wish to thank Fan Chung, Franklin Kenter and Choongbum Lee for fruitful discussions.


\bigskip

\begin{thebibliography}{99}

\bibitem{AKPSS} M. Ajtai, J. Koml\'{o}s, J. Pintz, J. Spencer and E. Szemer\'{e}di, Extremal uncrowded hypergraphs,
Journal of Combinatorial Theory A 32 (1982), no. 3, 321--335.

\bibitem{AKS} M. Ajtai, J. Koml\'{o}s,E. Szemer\'{e}di, A note on Ramsey numbers. J. Comb. Theory (Series A),
29 (1980), 354-–360.

\bibitem{AKSu} N. Alon, M. krivelevich, B. Sudakov, Coloring graphs with sparse neighborhoods, Journal
of Combinatorial Theory, Ser. B 77 (1999) 73--82.





\bibitem{BFM} T. Bohman, A. Frieze, D. Mubayi, Hypergraphs with independent neighborhoods, Combinatorica 30 (2010), no. 3, 277--293.

\bibitem{Cher} H. Chernoff, A Measure of Asymptotic Efficiency for Tests of a Hypothesis Based on the sum of Observations. Annals of Mathematical Statistics 23 (1952), no 4, 493–-507.

\bibitem{deBR} M. de Brandes, V. R\"{o}dl, Steiner triple systems with small maximal independent sets, Ars Combin. 17 (1984) 15--19.

\bibitem{DLR} R. Duke, H. Lefmann, V. R\"{o}dl, On uncrowded hypergraphs, Random Structures and Algorithms,
6, (1995) 209–-212.

\bibitem{E} P. Erd\H{o}s, Some unsolved problems in graph theory and combinatorial analysis in, in Combinatorial Mathematics and its Applications (Proc. Conf. Oxford, 1969), Academic Press, London, New York, 1971 pp. 97--109.



\bibitem{FM} A. Frieze, D. Mubayi, On the chromatic number of simple triangle-free triple systems,  Electron. J. Combin., 15(1) Research Paper 121, no. 27 (2008).

\bibitem{FM2} A. Frieze, D. Mubayi, Coloring simple hypergraphs,
submitted.

\bibitem{Fu} Z. F\"{u}redi, Maximal independent subsets in Steiner systems and in planar sets, SIAM Journal on Discrete Mathematics 4 (1991), 196--199.

\bibitem{FP} Z. F\"{u}redi, I.  Pal\'{a}sti, Arrangements of lines with a large number of triangles, Proceedings of the American Mathematical Society 92 no. 4 (1982) 561--566.

\bibitem{FMP} Z. F\"{u}redi, D. Mubayi, O. Pikhurko, Quadruple Systems with Independent Neighborhoods, J. Combin. Theory (Series A), 115 (2008) 1552--1560.

\bibitem{FPS} Z. F\"{u}redi, O. Pikhurko, M. Simonovits, On Triple Systems with Independent Neighborhoods, Combin. Probab. \& Comput., 14 (2005) 795--813.

\bibitem{K} J. H. Kim, The Ramsey Number $R(3,t)$ has order of magnitude $t^2/\log t$, Random Structures and Algorithms 7 (1995), 173--207

\bibitem{KPS} J. Koml\'{o}s, J. Pintz, E. Szemer\'{e}di, A lower bound for Heilbronn's problem, J.
London Math. Soc. 25 (1982), no. 2, 13--24.


\bibitem{McDiarmid} C. McDiarmid, Concentration, Probabilistic Methods for Algorithmic Discrete Mathematics (1998), 1--46.

\bibitem{PR} K. Phelps, V. R\"{o}dl, Steiner triple systems with minimum independence number, Ars Combin. 21 (1986), 167-–172.

\bibitem{R} V. R\"{o}dl, On a packing and covering problem, European J. Combin. 5 (1985), 69--78.

\bibitem{RS} V. R\"{o}dl, E. \v{S}inajov\'{a}, Note on independent sets in Steiner systems.
Proceedings of the Fifth International Seminar on Random Graphs and Probabilistic Methods in Combinatorics and Computer Science (Poznan, 1991).
Random Structures Algorithms 5 (1994), no. 1, 183–-190.

\bibitem{S} J. Shearer, A note on the independence number of triangle-free graphs. Discrete Mathematics 46 no 1. (1983) 83--87.

\end{thebibliography}
\end{document}